\documentclass{amsart}

\newtheorem{theorem}{Theorem}
\newtheorem{lemma}[theorem]{Lemma}

\theoremstyle{definition}
\newtheorem{definition}[theorem]{Definition}

\theoremstyle{remark}

\newcommand{\R}{{\mathbb R}}

\newcommand{\Z}{{\mathbb Z}}
\newcommand{\C}{{\mathbb C}}
\newcommand{\K}{{\mathbb K}}

\newcommand{\dv}[2]{\langle #1,#2\rangle}

\begin{document}

\title{On kernel theorems for (LF)-spaces}

\author{A.~G.~Smirnov}
\address{I.~E.~Tamm Theory Department, P.~N.~Lebedev
Physical Institute, Leninsky prospect 53, Moscow 119991,
Russia} \email{smirnov@lpi.ru}
\thanks{The research was supported by the grants
RFBR~05-01-01049, INTAS~03-51-6346, LSS-4401.2006.2, and
the grant MK-1315.2006.1 of the President of Russian
Federation.}

\date{}

\begin{abstract}
A convenient technique for proving kernel theorems for
(LF)-spaces (countable inductive limits of Fr\'echet
spaces)is developed. The proposed approach is based on
introducing a suitable modification of the functor of the
completed inductive topological tensor product. Using such
modified tensor products makes it possible to prove kernel
theorems without assuming the completeness of the
considered (LF)-spaces. The general construction is applied
to proving kernel theorems for a class of spaces of entire
analytic functions arising in nonlocal quantum field
theory.
\end{abstract}
\maketitle

\section{Introduction}

Let $X$ and $Y$ be sets and $F$, $G$, and $H$ be locally
convex spaces consisting of functions defined on $X$, $Y$,
and $X\times Y$ respectively. Suppose the function
$(x,y)\to f(x)g(y)$ belongs to $H$ for any $f\in F$ and
$g\in G$. In applications, it is often important to find
out whether the following statement holds:
\begin{itemize}
  \item[(K)] For any separately continuous bilinear mapping
  $\psi$ from $F\times G$ to a Hausdorff complete space
  $\tilde H$, there is a unique continuous linear mapping
  $\tilde \psi\colon H\to \tilde H$ such that $\psi=\tilde\psi
  \circ\Phi$, where the bilinear mapping $\Phi\colon
  F\times G\to H$ is defined by the relation
  $\Phi(f,g)(x,y) = f(x)g(y)$.
\end{itemize}
Statements of this type and their analogues for multilinear
mappings are known as kernel theorems. In this paper, we
propose a convenient technique for proving results of this
type in the case, where all considered spaces are
(LF)-spaces (i.e., countable inductive limits of Fr\'echet
spaces).

In the language of topological tensor products~\cite{Grot},
(K) means that $H$ can be identified with a dense subspace
of the completion $F\hat\otimes_i G$ of $F\otimes G$ with
respect to the inductive\footnote{Recall that the inductive
(projective) topology on $F\otimes G$ is the strongest
locally convex topology on $F\otimes G$ such that the
canonical bilinear mapping $(f,g)\to f\otimes g$ is
separately continuous (resp., continuous).} topology. In
particular, to prove~(K), it suffices to show that $H$ can
be identified with $F\hat\otimes_i G$. However, this
requires proving the completeness of $H$, which may present
difficulty for concrete functional (LF)-spaces (and is
actually unnecessary, as we shall see). To circumvent this
problem, we introduce a notion of the semi-completed tensor
product of (LF)-spaces (see~Definition~\ref{d1}). The
semi-completed tensor product $F\tilde \otimes \,G$ of
(LF)-spaces $F$ and $G$ is an (LF)-space which is
canonically identified with a dense subspace of
$F\hat\otimes_i G$, and, therefore, we can prove~(K) by
identifying $H$ with $F\tilde \otimes\, G$. Moreover,
semi-completed tensor products possess a natural
associativity property which turns out to be very useful
for the treatment of multilinear mappings.

The paper is organized as follows. In Sec.~\ref{s2}, we
give the definition of semi-completed tensor products and
describe their basic properties. In Sec.~\ref{s3}, we find
simple conditions ensuring the coincidence of $F\tilde
\otimes G$ with a given space $H$ in the case, where $F$,
$G$, and $H$ are functional (LF)-spaces. In Sec.~\ref{s4},
we apply the obtained results to proving kernel theorems
for some spaces of entire analytic functions arising in
nonlocal quantum field theory (see~\cite{SS2}).

In Secs.~\ref{s2} and~\ref{s3}, we assume that the
considered vector spaces are either real or complex. The
ground field ($\R$ or $\C$) is denoted by $\K$.

\section{Semi-completed tensor products of (LF)-spaces}
\label{s2}

A Hausdorff locally convex space $F$ is called an
(LF)-space if there exist a sequence $(F_i)$ of Fr\'echet
spaces and a sequence $p_i$ of continuous mappings from
$F_i$ to $F$ such that $F=\bigcup_i \mathrm{Im}\,p_i$ and
the topology of $F$ coincides with the inductive topology
with respect to the mappings $p_i$. The following important
result was proved by Grothendieck (see~Th\'eor\`eme~A
of~\cite{Grot}).

\begin{lemma}
 \label{l0}
Let $F$ be a Hausdorff locally convex space, $(F_i)$ be a
sequence of Fr\'echet spaces, and $(p_i)$ be a sequence of
continuous mappings from $F_i$ to $F$. Let $\tilde F$ be a
Fr\'echet space and $p\colon \tilde F\to F$ be a continuous
mapping such that $\mathrm{Im}\,p\subset \bigcup_i
\mathrm{Im}\,p_i$. Then there is an index $i$ such that
$\mathrm{Im}\,p\subset \mathrm{Im}\,p_i$ and if $p_i$ is
injective, then there is a continuous mapping $\tilde
p\colon \tilde F\to F_i$ such that $p=p_i\circ \tilde p$.
\end{lemma}

It follows from Lemma~\ref{l0} that the topology of an
(LF)-space $F$ coincides with the inductive topology with
respect to an arbitrary countable family of continuous
mappings of Fr\'echet spaces into $F$ provided that the
images of these mappings cover $F$.

Let $F$ be a Hausdorff locally convex space and $\tilde F$
be a Fr\'echet space. We say that $\tilde F$ is a Fr\'echet
subspace of $F$ (notation $\tilde F\prec F$) if $\tilde
F\subset F$ and the inclusion mapping $\tilde F\to F$ is
continuous. The relation $\prec$ is a partial order on the
set $\mathfrak F_F$ of all Fr\'echet subspaces of $F$. In
fact, $\mathfrak F_F$ is a lattice: the least upper bound
$F_1\vee F_2$ of $F_1,F_2\in \mathfrak F_F$ is the space
$F_1+F_2$ endowed with the inductive topology with respect
to the inclusion mappings $F_{1,2}\to F_1+F_2$; the
greatest lower bound $F_1\wedge F_2$ is the space $F_1\cap
F_2$ endowed with the projective topology with respect to
the inclusion mappings $F_1\cap F_2\to F_{1,2}$. In
particular, the set $\mathfrak F_F$ is directed. By
Lemma~\ref{l0}, a countable set $U\subset \mathfrak F_F$ is
cofinal in $\mathfrak F_F$ if and only if
$F=\bigcup_{\tilde F\in U} \tilde F$. If $F$ is an
(LF)-space, then $\mathfrak F_F$ contains a countable
cofinal subset. Indeed, let a sequence $(F_i)$ of Fr\'echet
spaces and a sequence $(p_i)$ of continuous linear mappings
from $F_i$ to $F$ be such that $F=\bigcup_i
\mathrm{Im}\,p_i$. Let $\tilde F_i$ be the space
$\mathrm{Im}\,p_i$ endowed with the inductive topology with
respect to the mapping $p_i$. Then $\tilde F_i\in \mathfrak
F_F$ for all $i$ and, therefore, $\tilde F_i$ form the
required subset.

For a locally convex space $F$, we denote by $\hat F$
the Hausdorff completion of $F$.
Given
locally convex spaces $F_1,\ldots,F_n$, the tensor product
$F_1\otimes\ldots\otimes F_n$ endowed by the inductive
(resp., projective) topology will be denoted by
$F_1\otimes_i\ldots\otimes_i F_n$ (resp., by
$F_1\otimes_\pi\ldots\otimes_\pi F_n$). The Hausdorff
completions of $F_1\otimes_i\ldots\otimes_i F_n$ and
$F_1\otimes_\pi\ldots\otimes_\pi F_n$ will be denoted by
$F_1\hat\otimes_i\ldots\hat\otimes_i F_n$ and
$F_1\hat\otimes_\pi\ldots\hat\otimes_\pi F_n$ respectively.
If $F_1,\ldots,F_n$ are Fr\'echet spaces, then the
inductive and projective topologies on
$F_1\otimes\ldots\otimes F_n$ coincide. In this case, we
shall omit the indices $i$ and $\pi$ in the notation for
tensor products.

\begin{definition}
\label{d1} Let $F_1,\ldots,F_n$ be (LF)-spaces. The
semi-completed tensor product
$F_1\tilde\otimes\ldots\tilde\otimes F_n$ of
$F_1,\ldots,F_n$ is defined to be the Hausdorff space
associated with the inductive limit
\[
\varinjlim_{(\tilde F_1,\ldots,\tilde F_n)\in \mathfrak
F_{F_1}\times\ldots\times\mathfrak F_{F_n}}\,\tilde
F_1\hat\otimes\ldots\hat\otimes\tilde F_n.
\]
\end{definition}

Each of the sets $\mathfrak F_{F_1},\ldots,\mathfrak
F_{F_n}$ contains a countable cofinal subset. It hence
follows that $\mathfrak F_{F_1}\times\ldots\times\mathfrak
F_{F_n}$ also contains a countable cofinal subset and,
therefore, $F_1\tilde\otimes\ldots\tilde\otimes F_n$ is an
(LF)-space. For $\tilde F=(\tilde F_1,\ldots,\tilde F_n)\in
\mathfrak F_{F_1}\times\ldots\times\mathfrak F_{F_n}$, let
$\iota_{\tilde F}$ be the natural continuous linear mapping
from $\tilde F_1\otimes\ldots\otimes \tilde F_n$ to
$F_1\otimes_i\ldots\otimes_i F_n$ determined by the
inclusion mappings $\tilde F_i\to F_i$, $i=1,\ldots,n$, and
$\rho_{\tilde F}$ be the canonical continuous linear
mapping from $\tilde F_1\hat\otimes\ldots\hat\otimes\tilde
F_n$ to $F_1\tilde\otimes\ldots\tilde\otimes F_n$.

\begin{lemma} \label{l0a}
Let $F_1,\ldots,F_n$ be (LF)-spaces. There is a unique
continuous linear mapping $j\colon
F_1\otimes_i\ldots\otimes_i F_n\to
F_1\tilde\otimes\ldots\tilde\otimes F_n$ such that
\[
j\iota_{\tilde F}= \rho_{\tilde F} \lambda_{\tilde F},
\]
for any $\tilde F=(\tilde F_1,\ldots,\tilde F_n)\in
\mathfrak F_{F_1}\times\ldots\times\mathfrak F_{F_n}$,
where $\lambda_{\tilde F}$ is the canonical mapping from
$\tilde F_1\otimes\ldots\otimes \tilde F_n$ to $\tilde
F_1\hat\otimes\ldots\hat\otimes \tilde F_n$. The mapping
$j$ is a topological isomorphism of
$F_1\otimes_i\ldots\otimes_i F_n$ onto a dense subspace of
$F_1\tilde\otimes\ldots\tilde\otimes F_n$. There is a
unique continuous linear mapping $k\colon
F_1\tilde\otimes\ldots\tilde\otimes F_n\to
F_1\hat\otimes_i\ldots\hat\otimes_i F_n$ such that $k\circ
j$ is the canonical mapping $F_1\otimes_i\ldots\otimes_i
F_n\to F_1\hat\otimes_i\ldots\hat\otimes_i F_n$. The
mapping $k$ is a topological isomorphism of
$F_1\tilde\otimes\ldots\tilde\otimes F_n$ onto a dense
subspace of $F_1\hat\otimes_i\ldots\hat\otimes_i F_n$. For
every $f\in F_1\tilde\otimes\ldots\tilde\otimes F_n$, there
is a bounded subset $B$ of $F_1\otimes_i\ldots\otimes_i
F_n$ such that $f\in \overline{j(B)}$, where the bar means
closure.
\end{lemma}

The mapping $j$ described in Lemma~\ref{l0a} will be called
the canonical mapping from $F_1\otimes_i\ldots\otimes_i
F_n$ to $F_1\tilde\otimes\ldots\tilde\otimes F_n$. For
$f_1\in F_1,\ldots, f_n\in F_n$, we set
\[
f_1\tilde\otimes \ldots \tilde\otimes f_n =
j(f_1\otimes\ldots\otimes f_n).
\]
The $n$-linear mapping $(f_1,\ldots, f_n)\to
f_1\tilde\otimes \ldots \tilde\otimes f_n$ will be called
the canonical $n$-linear mapping from
$F_1\times\ldots\times F_n$ to
$F_1\tilde\otimes\ldots\tilde\otimes F_n$.

\begin{definition}\label{d2}
Let $F_1,\ldots,F_n$ be (LF)-spaces, $H$ be a locally
convex space, and $\varphi$ be an $n$-linear mapping from
$F_1\times\ldots\times F_n$ to $H$. We say that $\varphi$
is (F)-continuous if for any $\tilde F_1\in\mathfrak
F_{F_1},\ldots,\tilde F_n\in\mathfrak F_{F_n}$, there are a
Fr\'echet space $\tilde H$, a continuous $n$-linear mapping
$\psi\colon \tilde F_1 \times\ldots\times \tilde
F_n\to\tilde H$, and a continuous linear mapping $l\colon
\tilde H\to H$ such that $l\circ\psi$ coincides with the
restriction of $\varphi$ to $\tilde F_1 \times\ldots\times
\tilde F_n$.
\end{definition}

Obviously, every (F)-continuous $n$-linear mapping is
separately continuous. If $H$ is Hausdorff and complete,
then every separately continuous $n$-linear mapping
$\varphi\colon F_1\times\ldots\times F_n\to H$ is
(F)-continuous. Indeed, for any $\tilde F_1\in\mathfrak
F_{F_1},\ldots,\tilde F_n\in\mathfrak F_{F_n}$, the
restriction of $\varphi$ to $\tilde F_1 \times\ldots\times
\tilde F_n$ can be decomposed as $l\circ\psi$, where $\psi$
is the canonical $n$-linear mapping from $\tilde F_1
\times\ldots\times \tilde F_n$ to $\tilde F_1
\hat\otimes\ldots\hat\otimes \tilde F_n$ and $l\colon
\tilde F_1 \hat\otimes\ldots\hat\otimes \tilde F_n\to H$ is
a continuous linear mapping.

It easily follows from the above definitions that the
canonical $n$-linear mapping from $F_1\times\ldots\times
F_n$ to $F_1\tilde\otimes\ldots\tilde\otimes F_n$ is
(F)-continuous. Moreover, the space
$F_1\tilde\otimes\ldots\tilde\otimes F_n$ has the following
universal property:

\begin{lemma}\label{l0b}
Let $F_1,\ldots,F_n$ be (LF)-spaces, $H$ be a Hausdorff
locally convex space, and $\varphi$ be an (F)-continuous
$n$-linear mapping from $F_1\times\ldots\times F_n$ to $H$.
Then there is a unique continuous linear mapping $l\colon
F_1\tilde\otimes\ldots\tilde\otimes F_n\to H$ such that
\[
\varphi(f_1,\ldots,f_n) =
l(f_1\tilde\otimes\ldots\tilde\otimes f_n),\quad f_1\in
F_1,\ldots, f_n\in F_n.
\]
\end{lemma}

\begin{theorem}\label{t0}
Let $F_1,\ldots,F_n$ be (LF)-spaces. For every $1\leq m<
n$, there is a unique topological isomorphism
$F_1\tilde\otimes\ldots\tilde\otimes F_n\simeq
(F_1\tilde\otimes\ldots\tilde\otimes F_m)\tilde\otimes\,
(F_{m+1}\tilde\otimes\ldots\tilde\otimes F_n)$ taking
$f_1\tilde\otimes\ldots\tilde\otimes f_n$ to
$(f_1\tilde\otimes\ldots\tilde\otimes f_m)\tilde\otimes\,
(f_{m+1}\tilde\otimes\ldots\tilde\otimes f_n)$ for any
$f_1\in F_1,\ldots, f_n\in F_n$.
\end{theorem}

\section{Tensor products of functional (LF)-spaces}
\label{s3}

Given a locally convex space $F$, we denote by $F'$ and
$\dv{\cdot}{\cdot}$ the continuous dual of $F$ and the
canonical bilinear form on $F'\times F$ respectively. We
denote by $F'_\sigma$, $F'_\tau$, and $F'_b$ the space $F'$
endowed with its weak topology, Mackey topology, and strong
topology respectively.

Let $X$ and $Y$ be sets and $F$ and $G$ be locally convex
spaces consisting of functions defined on $X$ and $Y$
respectively. We denote by $\Pi(F,G)$ the linear space
consisting of all functions $h(x,y)$ on $X\times Y$ such
that $h(x,\cdot)\in G$ for every $x\in X$ and the function
$h_v(x)=\dv{v}{h(x,\cdot)}$ belongs to $F$ for every $v\in
G'$.

\begin{lemma}
 \label{l1}
Let $X$ and $Y$ be sets and $F$, $G$, and $H$ be Hausdorff
complete locally convex spaces consisting of scalar functions
defined on $X$, $Y$, and $X\times Y$ respectively. Let $F$ be
B-complete, $G$ be nuclear, and the topologies of $F$, $G$, and
$H$ be stronger than that of simple convergence. Suppose the
function $(x,y)\to f(x)g(y)$ on $X\times Y$ belongs to $H$ for
every $f\in F$ and $g\in G$ and the bilinear mapping $\Phi\colon
F\times G\to H$ taking $(f,g)$ to this function is continuous.
Then $\Phi$ induces an injective continuous linear  mapping
$F\widehat{\otimes}_\pi G\to H$ whose image coincides with
$\Pi(F,G)$.
\end{lemma}

\begin{proof}
Without loss of generality, we can assume that $H=\K^{X\times Y}$.
Let $\Phi_*\colon F\widehat{\otimes}_\pi G \to H$ be the
continuous linear mapping determined by $\Phi$. As usual, let
$\mathfrak B_e(F'_\sigma, G'_\sigma)$ denote the space of
separately continuous bilinear forms on $F'_\sigma\times
G'_\sigma$ equipped with the biequicontinuous convergence topology
(i.e., the topology of the uniform convergence on the sets of the
form $A\times B$, where $A$ and $B$ are equicontinuous sets in
$F'$ and $G'$ respectively). Let $S$ be the natural continuous
linear mapping $F\widehat{\otimes}_\pi G \to \mathfrak
B_e(F'_\sigma, G'_\sigma)$ which takes $f\otimes g$ to the
bilinear form $(u,v)\to \dv{u}{f}\dv{v}{g}$. Since $G$ is nuclear,
$S$ is a topological isomorphism (see~\cite{Grot}, Chapitre~2,
Th\'eor\`eme~6 or~\cite{Schaefer}). Further, let $T$ be the linear
mapping $\mathfrak B_e(F'_\sigma, G'_\sigma)\to \K^{X\times Y}$
defined by the relation $(Tb)(x,y)=b(\delta_x,\delta_y)$, $b\in
\mathfrak B_e(F'_\sigma, G'_\sigma)$ (if $x\in X$ and $y\in Y$,
then $\delta_x$ and $\delta_y$ are the linear functionals on $F$
and $G$ such that $\dv{\delta_x}{f}=f(x)$ and
$\dv{\delta_y}{g}=g(y)$; they are continuous because the
topologies of $F$ and $G$ are stronger than the topology of simple
convergence). Obviously, $T$ is continuous and $TS$ coincides with
$\Phi_*$ on $F\otimes G$. By continuity, we have $\Phi_*=TS$
everywhere on $F\widehat{\otimes}_\pi G$. Moreover, $T$ is
injective because $\delta$-functionals are weakly dense in $F'$
and $G'$. To prove the statement, we therefore have to show that
$\mathrm{Im}\,T=\Pi(F,G)$.

Let $b\in \mathfrak B_e(F'_\sigma, G'_\sigma)$ and $h=Tb$. The
bilinear form $b$ determines two linear mappings $L_1\colon F'\to
G$ and $L_2\colon G'\to F$ such that $b(u,v)=\dv{v}{L_1
u}=\dv{u}{L_2 v}$ for any $u\in F'$ and $v\in G'$. For $x\in X$
and $y\in Y$, we have
$h(x,y)=\dv{\delta_y}{L_1\delta_x}=(L_1\delta_x)(y)$, i.e.,
$h(x,\cdot)=L_1\delta_x$. Further, for $v\in G'$ and $x\in X$, we
have $h_v(x)=\dv{v}{h(x,\cdot)}=
\dv{v}{L_1\delta_x}=\dv{\delta_x}{L_2v}=(L_2 v)(x)$. Hence
$h_v=L_2 v$ belongs to $G$ and, therefore, $h\in \Pi(F,G)$. Thus,
we have the inclusion $\mathrm{Im}\,T\subset\Pi(F,G)$.

We now prove the converse inclusion. Let $h\in \Pi(F,G)$ and
$L\colon G'_\tau\to F$ be the linear mapping taking $v\in G'$ to
$h_v$. We claim that the graph $\mathcal G$ of $L$ is closed. It
suffices to show that if an element of the form $(0,f)$ belongs to
the closure $\bar{\mathcal G}$ of $\mathcal G$, then $f=0$.
Suppose the contrary that there is $f_0\in F$ such that $f_0\ne 0$
and $(0,f_0)\in \bar{\mathcal G}$. Let $x_0\in X$ be such that
$f_0(x_0)\ne 0$ and let the neighborhood $U$ of $f_0$ be defined
by the relation $U=\{f\in F : |\dv{\delta_{x_0}}{f-f_0}|<
|f_0(x_0)|/2\}$. Let $V=\{v\in G' :
|\dv{v}{h(x_0,\cdot)}|<|f_0(x_0)|/2\}$. If $f\in U$ and $v\in V$,
then we have $|h_v(x_0)|<|f_0(x_0)|/2<|f(x_0)|$. Hence the
neighborhood $V\times U$ of $(0,f_0)$ does not intersect $\mathcal
G$. This contradicts to the assumption that $(0,f_0)\in
\bar{\mathcal G}$, and our claim is proved. Being nuclear and
complete, $G$ is semireflexive and hence $G'_\tau$ is barrelled.
We can therefore apply the closed graph theorem (\cite{Schaefer},
Theorem~IV.8.5) and conclude that $L$ is continuous. Let the
bilinear form $b$ on $F'\times G'$ be defined by the relation
$b(u,v)=\dv{u}{h_v}$. The continuity of $L$ implies that $b\in
\mathfrak B_e(F'_\sigma, G'_\sigma)$. Since $b(\delta_x,
\delta_y)=h(x,y)$, we have $h=Tb$. Thus, $h\in\mathrm{Im}\,T$ and
the lemma is proved.
\end{proof}

\begin{theorem}\label{t1}
Let $X$ and $Y$ be sets, $F$, $G$, and $H$ be (LF)-spaces
consisting of scalar functions on $X$, $Y$, and $X\times Y$
respectively. Suppose the topologies of $F$, $G$, and $H$
are stronger than that of simple convergence and $G$ can be
covered by a countable family of its nuclear Fr\'echet
subspaces. Let the following conditions be satisfied:
\begin{itemize}
\item[($i$)] For every $f\in F$ and $g\in G$, the function
$(x,y)\to f(x)g(y)$ on $X\times Y$ belongs to $H$ and the
bilinear mapping $\Phi\colon F\times G\to H$ taking $(f,g)$
to this function is (F)-continuous.
\item[($ii$)] For any $h\in H$, one can find Fr\'echet subspaces
$\tilde F\subset F$ and $\tilde G\subset G$ such that
$h(x,\cdot)\in \tilde G$ for every $x\in X$ and the
function $h_v(x)=\dv{v}{h(x,\cdot)}$ belongs to $\tilde F$
for every $v\in \tilde G'$.
\end{itemize}
Then $\Phi$ induces the topological isomorphism
$F\tilde\otimes \,G\simeq H$.
\end{theorem}

\begin{proof}
Let $\Phi_*\colon F\tilde{\otimes}\, G\to H$ be the
continuous linear mapping determined by $\Phi$. For $\tilde
F\in\mathfrak F_F$ and $\tilde G\in\mathfrak F_G$, let
$\Phi_*^{\tilde F,\tilde G}$ be the continuous linear
mapping $\tilde F\hat\otimes\, \tilde G\to H$ determined by
the restriction of $\Phi$ to $\tilde F\times \tilde G$. We
have
\begin{equation}\label{eq}
\Phi_*^{\tilde F,\tilde G}= \Phi_*\,\rho_{\tilde F,\tilde
G},
\end{equation}
where $\rho_{\tilde F,\tilde G}$ is the canonical mapping
from $\tilde F\hat\otimes\, \tilde G$ to
$F\tilde{\otimes}\, G$. Let $\mathcal N$ be the subset of
$\mathfrak F_G$ consisting of all nuclear Fr\'echet
subspaces of $G$. By the assumption, $\mathcal N$ is
cofinal in $\mathfrak F_G$. Let $h\in H$. Condition~(ii)
means that $h\in\Pi(\tilde F,\tilde G)$ for some $\tilde
F\in\mathfrak F_F$ and $\tilde G\in\mathfrak F_G$. Let
$\tilde G_1\in \mathcal N$ be such that $\tilde
G\prec\tilde G_1$. Then we have $\Pi(\tilde F,\tilde
G)\subset \Pi(\tilde F,\tilde G_1)$ and, therefore, $h\in
\Pi(\tilde F,\tilde G_1)$. By Lemma~\ref{l1}, we have
$\Pi(\tilde F,\tilde G_1)=\mathrm{Im}\, \Phi_*^{\tilde
F,\tilde G_1}$ and in view of~(\ref{eq}) we conclude that
$h\in \mathrm{Im}\,\Phi_*$. Thus, $\Phi_*$ is surjective.
By Lemma~\ref{l1}, the mappings $\Phi_*^{\tilde F,\tilde
G}$ are injective for any $\tilde F\in\mathfrak F_F$ and
$\tilde G\in \mathcal N$. Since $\mathcal N$ is cofinal in
$\mathfrak F_G$, it follows from~(\ref{eq}) that the
restriction of $\Phi_*$ to the image of $\rho_{\tilde
F,\tilde G}$ is injective for any $\tilde F\in\mathfrak
F_F$ and $\tilde G\in\mathfrak F_G$. This implies the
injectivity of $\Phi_*$ because the spaces
$\mathrm{Im}\,\rho_{\tilde F,\tilde G}$ cover $F\tilde
\otimes\, G$. We have thus proved that $\Phi_*$ is a
one-to-one mapping from $F\tilde \otimes\, G$ onto $H$.
Since both $F\tilde \otimes\, G$ and $H$ are (LF)-spaces,
we can apply the open mapping theorem (\cite{Grot},
Th\'eor\`eme~B) and conclude that $\Phi_*$ is a topological
isomorphism. The theorem is proved.
\end{proof}

\section{Applications to spaces of analytic functions}
\label{s4}

Let $U$ be a cone in $\R^k$. We say that a cone $W$ is a
conic neighborhood of $U$ if $W$ has an open
projection\footnote{The projection $\Pr W$ of a cone
$W\subset \R^k$ is by definition the intersection of $W$
with the unit sphere in $\R^k$; the projection of $W$ is
meant to be open in the topology of this sphere. } and
contains $U$.

\begin{definition} \label{d3}
Let $\beta(s)$ be continuous monotone indefinitely
increasing convex function on the semi-axis $s\geq 0$ and
$U$ be a nonempty cone in $\R^k$. For $B>0$, the Fr\'echet
space $\mathcal E^{\beta, B}(U)$ consists of entire
analytic functions on $\C^k$ having the finite norms
\[
\|f\|_{U,N,B'}=\sup_{z=x+iy\in\C^k}|f(z)|(1+|x|)^N
\exp\left(-\beta(B'|y|)-\beta(B'\delta_U(x))\right)
\]
for any $B'>B$ and any nonnegative integer $N$, where
$\delta_U(x)=\inf_{x'\in U}|x-x'|$ is the distance from $x$
to $U$. The space $\mathcal E^{\beta}(U)$ is defined by the
relation $\mathcal E^{\beta}(U)=\bigcup_{B>0,\,W\supset U}
\mathcal E^{\beta, B}(W)$, where $W$ runs over all conic
neighborhoods of $U$ and the union is endowed with the
inductive limit topology.
\end{definition}

The spaces $\mathcal E^\beta(\R^k)$ proved to be useful for
the analysis of infinite series in the Wick powers of free
fields converging to nonlocal fields~\cite{SS2}. If
$\beta(s)=s^{1/(\nu-1)}$, $\nu<1$, then $\mathcal
E^\beta(\R^k)$ coincides with the Gelfand-Shilov space
$S^\nu(\R^k)$ (see~\cite{GS} for the definition and
properties of Gelfand-Shilov spaces; to avoid confusion, we
use the notation $S^\nu$ instead of the standard
$S^\beta$). The spaces $\mathcal E^{\beta}(U)$ over cones
are introduced to describe the localization properties of
analytic functionals belonging to $\mathcal
E^{\prime\beta}(\R^k)$. More precisely, a closed cone
$K\subset \R^k$ is called a carrier cone of $u\in \mathcal
E^{\prime\beta}(\R^k)$ if $u$ has a continuous extension to
$\mathcal E^{\prime\beta}(K)$. The notion of carrier cone
replaces the notion of support of a generalized function
for elements of $E^{\prime\beta}(\R^k)$. In particular,
every $u\in \mathcal E^{\prime\beta}(\R^k)$ has a uniquely
determined minimal carrier cone (this was proved
in~\cite{Soloviev} for the case of the space $S^0$; the
general case can be treated in the same way using the
estimates for plurisubharmonic functions obtained
in~\cite{Smirnov}).

Here, we shall prove a kernel theorem for the spaces
$\mathcal E^{\beta}(U)$. For this, we introduce, in
addition to $\mathcal E^{\beta}(U)$, similar spaces
associated with finite families of cones.

\begin{definition} \label{d4}
Let $U_1,\ldots,U_n$ be nonempty cones in
$\R^{k_1},\ldots,\R^{k_n}$ respectively. We define the
space $\mathcal E^\beta(U_1,\ldots,U_n)$ by the relation
\[
\mathcal E^\beta(U_1,\ldots,U_n) =
\bigcup_{B>0,\,W_1\supset U_1,\ldots,W_n\supset U_n}
\mathcal E^{\beta,B}(W_1\times\ldots\times W_n),
\]
where the union is taken over all conic neighborhoods
$W_1,\ldots,W_n$ of $U_1,\ldots,U_n$ and is endowed with
the inductive limit topology.
\end{definition}

In what follows, we assume for definiteness that the norm
$|\cdot|$ on $\R^k$ is uniform: $|x|=\max_{1\leq j\leq
k}|x_j|$. For any cone $U$, the space
$\mathcal{E}^{\beta,B}(U)$ belongs to the class of the
spaces $\mathcal H(M)$ introduced in~\cite{SS1}. The space
$\mathcal H(M)$ is defined\footnote{The definition of
$\mathcal H(M)$ given here is slightly less general than
that in~\cite{SS1}, but it is sufficient for our purposes.}
by a family $M=\{M_\gamma\}_{\gamma\in\Gamma}$ of strictly
positive continuous functions on $\C^k$ and consists of all
entire analytic functions on $\C^k$ with the finite norms
\[
\sup_{z\in \C^k} M_\gamma(z)|f(z)|.
\]
It is supposed that (a) for every
$\gamma_1,\gamma_2\in\Gamma$, one can find
$\gamma\in\Gamma$ and $C>0$ such that $M_\gamma\geq C
(M_{\gamma_1}+M_{\gamma_2})$, and (b) there is a countable
set $\Gamma'\subset \Gamma$ with the property that for
every $\gamma\in \Gamma$, one can find $\gamma'\in \Gamma'$
and $C>0$ such that $C M_\gamma\leq M_{\gamma'}$ (a family
of functions on $\C^k$ satisfying~(a) and ~(b) will be
called a defining family of functions). Given a cone
$U\subset \R^k$ and $B>0$, we define the family $M^{U,B}$
of functions on $\C^k$ indexed by the set $\Z_+\times
(B,\infty)$ by the relation
\begin{equation}\label{xxx}
M^{U,B}_{N,B'}(x+iy)= (1+|x|)^{N}
\exp\left(-\beta(B'|y|)-\beta(B'\delta_U(x))\right),\quad
N\in\Z_+,\,B'>B.
\end{equation}
Then all above conditions are satisfied and we have
\begin{equation}\label{id}
\mathcal H(M^{U,B})= \mathcal{E}^{\beta,B}(U).
\end{equation}
By Lemma~12 of~\cite{SS1}, $\mathcal H(M)$ is a nuclear
Fr\'echet space if the following conditions are satisfied:
\begin{itemize}
\item[(I)] For every $\gamma\in \Gamma$, there
is $\gamma'\in \Gamma$ such that
$M_\gamma(z)/M_{\gamma'}(z)$ is integrable on $\C^k$ and
tends to zero as $|z|\to\infty$.
\item[(II)] For every $\gamma\in \Gamma$, there are
$\gamma'\in \Gamma$, a neighborhood of the origin $\mathcal
B$ in $\C^k$, and $C>0$ such that $M_\gamma(z)\leq
CM_{\gamma'}(z+\zeta)$ for any $z\in \C^k$ and $\zeta\in
\mathcal B$.
\end{itemize}
It is straightforward to verify that the family $M^{U,B}$
satisfies~(I) and~(II). The space
$\mathcal{E}^{\beta,B}(U)$ is therefore nuclear for any
$B>0$ and cone $U\subset\R^k$.

Let $M=\{M_\gamma\}_{\gamma\in \Gamma}$ and
$N=\{N_\omega\}_{\omega\in \Omega}$ be defining families of
functions on $\C^{k_1}$ and $\C^{k_2}$ respectively. We
denote by $M\otimes N$ the family formed by the functions
\[
(M\otimes N)_{\gamma\omega}(z_1,z_2)=M_\gamma(z_1)
N_\omega(z_2),\quad (\gamma,\omega)\in \Gamma\times \Omega.
\]
Clearly, if $M\otimes N$ is a defining family of functions.
The following result was proved in~\cite{SS1}.

\begin{lemma}
 \label{l3}
Let $M=\{M_\gamma\}_{\gamma\in \Gamma}$ and
$N=\{N_\omega\}_{\omega\in \Omega}$ be defining families of
functions on $\C^{k_1}$ and $\C^{k_2}$ respectively and let
$h\in \mathcal H(M\otimes N)$. Suppose $N$ satisfies
conditions~{\rm (I)} and~{\rm (II)}. Then $h(z,\cdot)\in
\mathcal H(N)$ for every $z\in \C^{k_1}$ and the function
$h_v(z)=\dv{v}{h(z,\cdot)}$ belongs to $\mathcal H(M)$ for
all $v\in \mathcal H'(N)$.
\end{lemma}

Let $V_1\subset\R^{k_1}$ and $V_1\subset\R^{k_1}$ be
nonempty cones. It follows from~(\ref{xxx}) and the
monotonicity of $\beta$ that
\begin{equation}\nonumber
M^{V_1\times V_2,2B}_{N,2B'}(z_1,z_2)\leq
M^{V_1,B}_{N,B'}(z_1)\,M^{V_2,B}_{N,B'}(z_2) \leq
M^{V_1\times V_2,B}_{2N,B'}(z_1,z_2),\quad
z_{1,2}\in\C^{k_{1,2}},
\end{equation}
for any $B'>B>0$ and $N\in\Z_+$. We hence have continuous
inclusions
\begin{equation}\label{zzz}
\mathcal E^{\beta,B}(V_1\times V_2)\subset \mathcal
H(M^{V_1,B}\otimes M^{V_2,B})\subset \mathcal
E^{\beta,2B}(V_1\times V_2).
\end{equation}

\begin{lemma}\label{l4}
Let $U_1,\ldots,U_n$ be nonempty cones in
$\R^{k_1},\ldots,\R^{k_n}$ respectively. Let $1\leq m<n$
and $\Phi\colon \mathcal E^\beta(U_1,\ldots,U_m)\times
\mathcal E^\beta(U_{m+1},\ldots,U_n)\to \mathcal
E^\beta(U_1,\ldots,U_n)$ be the bilinear mapping defined by
the relation
\[
\Phi(f,g)(z_1,\ldots,z_n) =
f(z_1,\ldots,z_m)\,g(z_{m+1},\ldots,z_n).
\]
Then $\Phi$ is (F)-continuous and induces the topological
isomorphism
\[
\mathcal E^\beta(U_1,\ldots,U_m)\tilde\otimes \,\mathcal
E^\beta(U_{m+1},\ldots,U_n)\simeq \mathcal
E^\beta(U_1,\ldots,U_n).
\]
\end{lemma}
\begin{proof}
We check that the spaces $F = \mathcal
E^\beta(U_1,\ldots,U_m)$, $G = \mathcal
E^\beta(U_{m+1},\ldots,U_n)$, and $H = \mathcal
E^\beta(U_1,\ldots,U_n)$ and the bilinear mapping $\Phi$
satisfy conditions~(i) and~(ii) of Theorem~\ref{t1}. Let
$B>0$ and $W_1,\ldots,W_n$ be conic neighborhoods of
$U_1,\ldots,U_n$ respectively. Let $V_1 =
W_1\times\ldots\times W_m$ and $V_2 =
W_{m+1}\times\ldots\times W_n$. In view of~(\ref{id}) it
follows from~(\ref{zzz}) that $\Phi$ induces a continuous
bilinear mapping $\mathcal E^{\beta,B}(V_1)\times \mathcal
E^{\beta,B}(V_2)\to \mathcal E^{\beta,B}(V_1\times V_2)$.
This implies that $\Phi$ is (F)-continuous and~(i) is
fulfilled. Let $h\in \mathcal
E^{\beta,B}(W_1\times\ldots\times W_n)$. It follows
from~(\ref{zzz}) and Lemma~\ref{l3} that (ii) will be
satisfied if we set $\tilde F = \mathcal E^{\beta,B}(V_1)$
and $\tilde G = \mathcal E^{\beta,B}(V_2)$. The lemma is
proved.
\end{proof}

The next result follows from Lemma~\ref{l4} and
Theorem~\ref{t0} by induction on $n$.

\begin{theorem}\label{t2}
Let $U_1,\ldots,U_n$ be nonempty cones in
$\R^{k_1},\ldots,\R^{k_n}$ respectively and $\Phi\colon
\mathcal E^\beta(U_1)\times\ldots\times \mathcal
E^\beta(U_n)\to \mathcal E^\beta(U_1,\ldots,U_n)$ be the
bilinear mapping defined by the relation
\[
\Phi(f_1,\ldots,f_n)(z_1,\ldots,z_n) = f_1(z_1)\ldots
f_n(z_n).
\]
Then $\Phi$ is (F)-continuous and induces the topological
isomorphism
\[
\mathcal E^\beta(U_1)\tilde\otimes\ldots\tilde\otimes
\,\mathcal E^\beta(U_n)\simeq \mathcal
E^\beta(U_1,\ldots,U_n).
\]

\end{theorem}

\end{document}